\documentclass[11pt]{amsart}
\usepackage{natbib}

\newtheorem{theorem}{Theorem}
\newtheorem{lemma}{Lemma}
\newtheorem{corollary}{Corollary}
\theoremstyle{remark}
\newtheorem{example}{Example}

\newcommand{\nor}{\triangleleft} 
\newcommand{\cF}{\mathcal{F}}
\newcommand{\cV}{\mathcal{V}}

\begin{document}

\title[Loops and the Lagrange property]{Loops and the Lagrange property}

\author[O.~Chein]{Orin Chein}
\address{Department of Mathematics,
Temple University,
1805 N. Broad St.,
Philadelphia, PA 19122 USA}
\email{orin@math.temple.edu}
\urladdr{http://math.temple.edu/\symbol{126}orin}

\author[M.~K.~Kinyon]{Michael K. Kinyon}
\address{Department of Mathematics,
Western Michigan University,
Kalamazoo, MI~49008-5248 USA}
\email{mkinyon@wmich.edu}
\urladdr{http://unix.cc.wmich.edu/\symbol{126}mkinyon}

\author[A.~Rajah]{\\ Andrew Rajah}
\address{School of Mathematical Sciences,
Universiti Sains Malaysia,
11800 USM Penang, Malaysia}
\email{andy@cs.usm.my}
\urladdr{http://www.mat.usm.my/math/staff2.htm{\#}AR}

\author[P.~Vojt\v{e}chovsk\'y]{Petr Vojt\v{e}chovsk\'y}
\address{Department of Mathematics,
Iowa State University,
Ames, IA~50011 USA}
\email{petr@iastate.edu}
\urladdr{http://www.vojtechovsky.com}


\begin{abstract}
Let $\cF$ be a family of finite loops closed under subloops and factor loops.
Then every loop in $\cF$ has the strong Lagrange property if and only if every
simple loop in $\cF$ has the weak Lagrange property. We exhibit several such
families, and indicate how the Lagrange property enters into the problem of
existence of finite simple loops.
\end{abstract}

\maketitle
\thispagestyle{empty}

\noindent
The two most important open problems in loop theory, namely the existence of a
finite simple Bol loop and the Lagrange property for Moufang loops, have been
around for more than 40 years. While we certainly have not solved these
problems, we show that they are closely related. Some of the ideas developed
here have been present in the loop-theoretical community, however, in a rather
vague form. We thus felt the need to express them more precisely and in a more
definite way.

We assume only basic familiarity with loops, not reaching beyond the
introductory chapters of \cite{pf}. All loops mentioned below are finite.

We begin with the crucial notion: the Lagrange property. A loop $L$ is said to
have the \emph{weak Lagrange property} if, for each subloop $K$ of $L$, $|K|$
divides $|L|$. It has the \emph{strong Lagrange property} if every subloop $K$
of $L$ has the weak Lagrange property.

A loop may have the weak Lagrange property but not the strong Lagrange
property. Four of the six nonisomorphic loops of order $5$ have elements of
order $2$ and hence fail to satisfy the weak Lagrange property. Let $K$ be one
of these loops. As noted in \cite[p.\ 13]{pf}, if $L$ is a loop of order $10$
having $K$ as a subloop and satisfying the property that every proper subloop
of $L$ is contained in $K$, then $L$ will have the weak but not the strong
Lagrange property. It is not difficult to construct a multiplication table for
such a loop.

Our main result depends on the following lemma, which is a restatement of
\cite[Lemma V.2.1]{br}.

\begin{lemma}\label{lem:bruck}
Let $L$ be a loop with a normal subloop $N$ such that
\begin{enumerate}
\item[(i)] $N$ has the weak $($resp.\ strong$)$ Lagrange property,
and
\item[(ii)] $L/N$ has the weak $($resp.\ strong$)$ Lagrange property.
\end{enumerate}
Then L has the weak $($resp.\ strong$)$ Lagrange property.
\end{lemma}

There are some classes of loops studied in the literature to which the lemma
applies directly. For each of these, the normal subloop in question is
associative. For instance, an easy induction shows that any solvable loop
satisfies the strong Lagrange property.  In particular, any (centrally)
nilpotent loop satisfies this property as well.

Before we turn to more specific examples, we recall a few definitions. Let $L$
be a loop and $x\in L$. When $x$ has a two-sided inverse, we denote it by
$x^{-1}$. A loop $L$ has the \emph{automorphic inverse property} if
$x^{-1}y^{-1}=(xy)^{-1}$ holds for every $x$, $y\in L$. A loop $L$ is
\emph{$($right$)$ Bol} (resp.\ \emph{left Bol}), if $((xy)z)y=x((yz)y)$ (resp.
$(x(yx))z=x(y(xz))$) holds for every $x, y, z\in L$. \emph{Moufang loops} are
loops that are both right Bol and left Bol. Since the concepts of right Bol
loop and left Bol loop are anti-isomorphic to each other, right Bol and left
Bol loops share the same algebraic properties. Thus, in what follows, when we
refer to Bol loops, the reader may think of left Bol or right Bol as he or she
sees fit. A Bol loop which has the automorphic inverse property is called a
\emph{Bruck loop}. Bruck loops of odd order are called \emph{B-loops}\cite[p.\
376]{gl1}. An \emph{A-loop} is a loop all of whose inner mappings are
automorphisms \cite{bp}. Finally, an \emph{$M_k$ loop} is a Moufang loop $L$
for which $L/\mathrm{Nuc}(L)$ has exponent $k-1$ and no smaller exponent, where
$\mathrm{Nuc}(L)$ is the nucleus of $L$ \cite{pf1, cp}.

\begin{example}
Let $L$ be a Moufang loop with an associative normal subloop $K$ such that
$L/K$ satisfies the strong Lagrange property. By the lemma, $L$ has the strong
Lagrange property.

As an example, if $L$ is an $M_k$ loop where $k = 2^m + 1$ and if $K =
\mathrm{Nuc}(L)$, then $L/K$ is a Moufang $2$-loop which, by \cite{gw} is
centrally nilpotent (cf. \cite{c1}).

As another example of this type, let $L$ be a Moufang loop with an associative
normal subloop $K$ such that $|L/K|$ is odd. In this case $L/K$ has the strong
Lagrange property by \cite[Thm.\ 2]{gl2}. For instance, if $L$ is an $M_k$ loop
with $k$ even, then $K = \mathrm{Nuc}(L)$ is an associative normal subloop such
that $L/K$ is of odd exponent, $k-1$, and so $L/K$ must have odd order (cf.
\cite{c1}).

Since every Moufang $A$-loop is an $M_4$ loop \cite[Cor.\ 2]{kkp}, and every
commutative Moufang loop is a Moufang $A$-loop \cite{bp}, it follows that every
Moufang $A$-loop, and in particular, every commutative Moufang loop has the
strong Lagrange property.
\end{example}

That every commutative Moufang loop has the strong Lagrange property is, in
fact, a well-known folk result, and follows from the central nilpotence of
these loops and the lemma.

\begin{example}
Let $L$ be a Bruck loop with an associative normal subloop $K$ such that
$|L/K|$ is odd. Since $L/K$ is a B-loop, it follows from \cite[Cor.\ 4]{gl1}
that $L/K$ has the strong Lagrange property. By the lemma, so does $L$. For
instance, since a commutative Moufang loop of odd order is obviously a B-loop,
this gives an alternative proof that every such loop has the strong Lagrange
property.
\end{example}

\begin{example}
The lemma also applies directly to those loops $L$ for which the derived
subloop $L'$ (i.e., the smallest subloop $L'$ such that $L/L'$ is an abelian
group) has the strong Lagrange property. For instance, let $L$ be a ``central''
Bol loop in the terminology of Kreuzer \cite{kr}, i.e., a Bol loop $L$ such
that $L'$ is contained in the center. These are just centrally nilpotent Bol
loops of nilpotence class 2, and thus these loops have the strong Lagrange
property. As another example, Bruck and Paige \cite{bp} showed that an A-loop
$L$ has the property that all of its loop isotopes are A-loops if and only if
$L/\mathrm{Nuc}(L)$ is a group, in other words, if and only if $L$ is nuclearly
nilpotent of class $2$. By the lemma, such an $L$ has the strong Lagrange
property.
\end{example}

We now come to our main result---the connection between simple loops and loops
satisfying the Lagrange property.

\begin{theorem}\label{thm:main}
Let $\cF$ be a nonempty family of finite loops such that
\begin{enumerate}
\item If $L \in \cF$ and $N \nor L$, then $N \in \cF$;
\item If $L \in \cF$ and $N \nor L$, then $L/N \in \cF$;
\item Every simple loop in $\cF$ has the weak Lagrange property.
\end{enumerate}
Then every loop in $\cF$ has the weak Lagrange property.
\end{theorem}

\begin{proof}
We proceed by induction on the order of loops in $\cF$. Note that (1) implies
that $\cF$ contains the trivial loop $\langle 1 \rangle$, for which the desired
conclusion is trivial. Now fix $L \in \cF$ and assume that the result holds for
all loops in $\cF$ of order less than $|L|$. If $L$ is simple, we are finished
by (3). Thus assume that $L$ is not simple, so that $L$ has a nontrivial proper
normal subloop $N$. Since $|N| < |L|$, (1) and the induction hypothesis imply
that $N$ has the weak Lagrange property. Since $|L/N| < |L|$, (2) and the
induction hypothesis imply that $L/N$ has the weak Lagrange property. By Lemma
\ref{lem:bruck}, $L$ has the weak Lagrange property.
\end{proof}

\begin{corollary}\label{coro:1}
Let $\cF$ be a nonempty family of finite loops such that
\begin{enumerate}
\item[(1')] If $L \in \cF$ and $K \leq L$, then $K \in \cF$;
\item[(2)] If $L \in \cF$ and $N \nor L$, then $L/N \in \cF$;
\item[(3)] Every simple loop in $\cF$ has the weak Lagrange property.
\end{enumerate}
Then every loop in $\cF$ has the strong Lagrange property.
\end{corollary}

\begin{proof}
Since (1') implies (1), Theorem \ref{thm:main} implies that every loop in $\cF$
has the weak Lagrange property. But then (1') yields the desired result.
\end{proof}

\begin{corollary}\label{coro:2}
Let $\cV$ be a variety of loops such that every simple loop in $\cV$ has the
weak Lagrange property. Then every loop in $\cV$ has the strong Lagrange
property.
\end{corollary}

Corollary \ref{coro:2} is of particular interest for those varieties of loops
for which there exists a classification of all simple loops. A prominent
example is the variety of Moufang loops, where it is known that every simple
nonassociative Moufang loop is isomorphic to a Paige loop (cf.\ \cite{paige},
\cite{liebeck}). It follows from Corollary \ref{coro:2} that if the weak
Lagrange property can be established for each of the Paige loops, then every
Moufang loop will have the strong Lagrange property.

There is one Paige loop for every finite field $GF(q)$; its order is
$q^3(q^4-1)$ when $q$ is even, and $q^3(q^4-1)/2$ when $q$ is odd \cite{paige}.
The weak Lagrange property for the smallest 120-element Paige loop has been
established in \cite{mgpm} and \cite{vo}. The next smallest Paige loop has
order $1080$. Thus based on published literature, we can state this result.

\begin{corollary}\label{coro:3}
Every Moufang loop of order less than $1080$ has the strong Lagrange property.
\end{corollary}

\begin{proof}
A simple Moufang loop of order less than $1080$ is a group or the smallest
Paige loop. The result follows from Corollary \ref{coro:1}.
\end{proof}

Incidentally, since no Paige loop is commutative, it follows from Corollary
\ref{coro:2} that, once again, every commutative Moufang loop has the strong
Lagrange property.

The authors have been informed by G.~E.~Moorehouse \cite{em} that by a computer
search, he has found that the Paige loop of order $1080$ satisfies the weak
Lagrange property. If we assume this to be correct, then we may state the
following.

\begin{corollary}\label{coro:4}
Every Moufang loop of order less than $16320$ has the strong Lagrange property.
\end{corollary}

Altogether, we have demonstrated that a loop has the strong Lagrange property
whenever it belongs to one of the following classes: Moufang loops with an
associative normal subloop of odd index, Bruck loops with an associative normal
subloop of odd index, loops whose derived subloops have the strong Lagrange
property, and Moufang loops of order less than 1080.

We conclude this paper with a couple of remarks on a potential application for
Corollary \ref{coro:2} to the existence of finite simple non-Moufang Bol loops.
This can be split into two problems: the existence of a finite simple Bruck
loop and the existence of a finite simple proper (non-Moufang, non-Bruck) Bol
loop.

First, to establish the existence of a finite simple Bruck loop, it would be
sufficient to find a Bruck loop which violated the weak Lagrange property, for
then Corollary \ref{coro:2} would imply the existence of a simple Bruck loop
which is not a cyclic group.

On the other hand, Corollary 2 might apply to the problem of the existence of
finite simple proper Bol loops if the weak Lagrange property is established for
all Paige loops. It would then be sufficient to find a Bol loop which violated
the weak Lagrange property, for by Corollary 2, there would exist a simple Bol
loop which violated that property. If all simple Moufang loops have the weak
Lagrange property, then the simple Bol loop in question will not be Moufang.

\bibliographystyle{plain}

\begin{thebibliography}{99}

\bibitem[Bruck(1958)]{br}
R.~H.~Bruck,
\textit{A Survey of Binary Systems},
Springer-Verlag, 1958; third printing, 1971.

\bibitem[Bruck and Paige(1956)]{bp}
R.~H.~Bruck and L.~J.~Paige,
Loops whose inner mappings are automorphisms,
\textit{Ann. of Math.} (2) \textbf{63} (1956) 308-323.

\bibitem[Chein(1973)]{c1}
O.~Chein,
Lagrange's theorem for $M_k$-loops,
\textit{Arch. Math. (Basel)} \textbf{24} (1973), 121-122.

\bibitem[Chein and Pflugfelder(1971)]{cp}
O.~Chein and H.~O.~Pflugfelder, On maps $x\rightarrow x\sp{n}$ and the
isotopy-isomorphy property of Moufang loops, \textit{Aequationes Math.} 6
(1971) 157-161.

\bibitem[Chein, Pflugfelder and Smith(1990)]{cps}
O.~Chein, H.~O.~Pflugfelder, and J.~D.~H.~Smith,
\textit{Quasigroups and Loops: Theory and Applications},
Sigma Series in Pure Math. \textbf{9}, Heldermann Verlag, 1990.

\bibitem[Glauberman(1964)]{gl1}
G.~Glauberman,
On loops of odd order I.
\textit{J. Algebra} \textbf{1} (1964) 374--396.

\bibitem[Glauberman(1968)]{gl2}
G.~Glauberman,
On loops of odd order II.
\textit{J. Algebra} \textbf{8} (1968) 393--414.

\bibitem[Glauberman and Wright(1968)]{gw}
G.~Glauberman and C.~R.~B.~Wright, Nilpotence of finite Moufang $2$-loops,
\textit{J. Algebra} 8 (1968) 415-417.

\bibitem[Guiliani and Milies(2000)]{mgpm}
M.~L.~Merlini Guiliani, C\'esar Polcino Milies, On the structure of the simple
Moufang loop $GLL(F_2)$, in \textit{Nonassociative algebra and its
applications}, the Fourth International Conference, Lecture Notes in Pure and
Applied Mathematics \textbf{211}, ed. R. Costa, A. Grishkov, H. Guzzo,~Jr., L.
A.~Peresi. Marcel Dekker, New York, 2000.

\bibitem[Kinyon, Kunen and Phillips(2002)]{kkp}
M.~Kinyon, K.~Kunen, and J.~D.~Phillips,
Every diassociative A-loop is Moufang,
\textit{Proc. Amer. Math. Soc.} \textbf{130} (2002) 619-624.

\bibitem[Kreuzer(1997)]{kr}
A.~Kreuzer,
Central Bol loops,
in \textit{Nearrings, nearfields and $K$-loops} (Hamburg,
1995), 301--310, Math. Appl. \textbf{426}, Kluwer Acad. Publ.,
Dordrecht, 1997.

\bibitem[Liebeck(1987)]{liebeck}
M.~W.~Liebeck, The classification of finite simple Moufang loops, \textit{Math.
Proc. Cambridge Philos. Soc.} \textbf{102} (1987) 33--47.

\bibitem[Moorehouse(2002)]{em}
G.~E.~Moorehouse, private communication.

\bibitem[Paige(1956)]{paige}
L.~J.~Paige,
A Class of Simple Moufang Loops,
\textit{Proc. Amer. Math. Soc.} \textbf{7} (1956) 471--482.

\bibitem[Pflugfelder(1970)]{pf1}
H.~Orlik-Pflugfelder,
A special class of Moufang loops,
\textit{Proc. Amer. Math. Soc.} \textbf{26} (1970) 583-586

\bibitem[Pflugfelder(1990)]{pf}
H.~O.~Pflugfelder,
{\it Quasigroups and Loops: Introduction},
Sigma Series in Pure Math. \textbf{8}, Heldermann Verlag,
Berlin, 1990.

\bibitem[Vojt\v{e}chovsk\'y(2001)]{vo}
P.~Vojt\v{e}chovsk\'y, \textit{Finite simple Moufang loops}, Ph.D.
Dissertation, Iowa State University, 2001. Available online at
www.vojtechovsky.com.

\end{thebibliography}

\end{document}